\documentclass[a4paper]{amsart}

\usepackage{amssymb}
\usepackage{amsthm}
\usepackage{amscd}
\usepackage{amsmath}
\usepackage[all]{xy}
\usepackage{graphicx}
\newtheorem{theorem}{Theorem}[subsection]
\newtheorem{claim}[theorem]{Claim}
\newtheorem{lemma}[theorem]{Lemma}
\newtheorem{proposition}[theorem]{Proposition}
\newtheorem{cor}[theorem]{Corollary}

\theoremstyle{definition}
\newtheorem{definition}[theorem]{Definition}
\newtheorem{remark}[theorem]{Remark}

\newcommand{\Z}{{\mathbb Z}}
\newcommand{\C}{{\mathbb C}}

\newcommand{\R}{{\mathbb R}}

\newcommand{\U}{{\rm U}}
\newcommand{\F}{{\mathbb F}}

\newcommand{\SO}{{\rm SO}}
\newcommand{\SU}{{\rm SU}}
\newcommand{\SL}{{\rm SL}}
\newcommand{\GL}{{\rm GL}}
\newcommand{\T}{\mathcal{T}}

\newcommand{\gtilde}{\widetilde{\mathfrak{g}}}
\newcommand{\grho}{\mathfrak{g}_\rho}
\newcommand{\mathfrakg}{\mathfrak{g}}

\newcommand{\trace}{{\rm Tr}\,}
\newfont{\bg}{cmr10 scaled\magstep4}

\begin{document}

\thanks{This research is partially supported by 
the 21st century COE program at Graduate School of Mathematical
Sciences, the University of Tokyo}

\title{
A relationship between 
the non-acyclic Reidemeister torsion and
a zero of the acyclic Reidemeister torsion
}


\author{Yoshikazu Yamaguchi}
\address{Graduate School of Mathematical Sciences, 
University of Tokyo, 3-8-1 Komaba Meguro, 
Tokyo 153-8914, Japan}
\email{shouji@ms.u-tokyo.ac.jp}

\date{}

\begin{abstract}
We show a relationship between 
the non-acyclic Reidemeister torsion
and
a zero of the acyclic Reidemeister torsion
for a $\lambda$-regular $\SU(2)$ or $\SL(2, \C)$-representation of 
a knot group.
Then we give a method to 
calculate the non-acyclic Reidemeister torsion of 
a knot exterior.
We calculate a new example and 
investigate the behavior of 
the non-acyclic Reidemeister torsion 
associated to a $2$-bridge knot 
and $\SU(2)$-representations of its knot group.
\end{abstract}

\maketitle



\section{Introduction.}
The Reidemeister torsion is an invariant 
for a CW-complex and a representation of its fundamental group.
In other words, this invariant associates with the local system
for a representation of the fundamental group.
Originally the Reidemeister torsion is defined 
if the local system is ${\it acyclic\/}$, i.e., 
all homology groups vanish.
However 
we can extend the definition of the Reidemeister torsion to non-acyclic cases 
\cite{Milnor66, Turaev1}.
In this paper, we focus on the non-acyclic cases.

It is known that the Fox calculus plays important roles in the
study of the Reidemeister torsion
\cite{dubois1, KL, Kitano, MilnorCov, Porti, Turaev1}.
The many results were obtained by using the Fox calculus
for the acyclic Reidemeister torsion.
In particular, there are important results related to the Alexander polynomial
in the knot theory \cite{KL, Kitano, MilnorCov, Turaev1}.
The Fox calculus is also important for non-acyclic cases \cite{dubois1, Porti}.
It is related to the cohomology theory of groups.

This paper contributes to the study of the non-acyclic Reidemeister torsion
by using the Fox calculus.
Our purpose is to apply the Fox calculus for the acyclic cases to
the study of the non-acyclic Reidemeister torsion
by using a relationship between the acyclic Reidemeister torsion 
and the non-acyclic one.
Our main theorem says that
the non-acyclic Reidemeister torsion for a knot exterior
is given by the differential coefficients of 
the twisted Alexander invariant of the knot.
The twisted Alexander invariant of a knot 
is the acyclic Reidemeister torsion and 
expressed as a one variable rational function \cite{Kitano}.
A conjecture due to J. Dubois and R. Kashaev \cite{DuboisKashaev} will be solved
in \cite{yamaguchi} by using our main theorem.

In the latter of this paper, we apply this relationship 
to study the Reidemeister torsion 
for the pair of a $2$-bridge knot and $\SU(2)$-representation of its knot group.
We give an explicit expression of the non-acyclic Reidemeister torsion 
associated to $5_2$ knot.
This is a new example of calculation of the non-acyclic Reidemeister torsion.
Furthermore, 
we investigate where the non-acyclic Reidemeister torsion associated to 
a $2$-bridge knot has critical points.
Note that
the non-acyclic Reidemeister torsion is parametrized by the representations of 
a knot group.
Moreover 
this Reidemeister torsion turns into 
a function on the character variety of the knot group. 
We will see that 
the critical points of the non-acyclic Reidemeister torsion associated to
a $2$-bridge knot are binary dihedral representations and 
these representations are related to 
the geometry of the character variety of a $2$-bridge knot group.

\medskip

This paper is organized as follows.
In Section \ref{Review_twisted_torsion}, we review the Reidemeister torsion.
In particular, we give the notion of the non-acyclic Reidemeister torsion of 
knot exteriors \cite{dubois1,Porti}.

Section \ref{Main_theorem} includes our main theorem on
a relationship between the non-acyclic Reidemeister torsion 
and the twisted Alexander invariant for knot exteriors.
We give a formula of the non-acyclic Reidemeister torsion for a knot exterior
by using a Wirtinger presentation of a knot group.

In Section \ref{applications}, 
we apply the results of Section \ref{Main_theorem}
to study the non-acyclic Reidemeister torsion for a $2$-bridge knot group and 
$\SU(2)$-representation of its knot group.

\section{Review on the non-abelian twisted Reidemeister torsion}
\label{Review_twisted_torsion}
\subsection{Notation}
In this paper, we use the following notations.
\begin{itemize}
\item
$\F$ is the field $\R$ or $\C$. 
\item
$G$ is the Lie group $\SU(2)$ (resp. $\SL(2, \C))$ 
if $\F$ is $\R$ (resp. $\C$).  
The symbol $\mathfrakg$ denotes the Lie algebra of $G$.

\item
$Ad$ denotes the adjoint action of $G$ to the Lie group $\mathfrakg$.
\item
$(\,  , \, )_{\mathfrakg} : \mathfrakg \times \mathfrakg \to \F$ is 
a product on the $\mathfrakg$, 
which is defined by 
$
(X, Y)_{\mathfrakg} = {\trace ( ^t \! X \bar Y)}.
$
\item
$V$ denotes an $n$-dimensional vector space over $\F$.
\item
For two ordered bases $\bf a$ and $\bf b$ in a vector space, 
we denote by $({\bf a}/{\bf b})$ the base-change matrix 
from ${\bf b}$ to ${\bf a}$ satisfying ${\bf a}={\bf b}({\bf a}/{\bf b})$.
We write simply $[{\bf a}/{\bf b}]$ for 
the determinant $\det ({\bf a}/{\bf b})$ of $({\bf a}/{\bf b})$.
We deal with ordered bases in this paper.
\end{itemize}
\subsection{Torsion of a chain complex} 
We recall the definition of the torsion.

Let 
$C_* = 
(0 \to C_n 
   \xrightarrow{\partial_n} C_{n-1} 
   \xrightarrow{\partial_{n-1}} 
   \cdots 
   \xrightarrow{\partial_1} C_0 
   \to 0)$
be a chain complex over $\F$.
For each $i$ let $Z_i$ denote the kernel of $\partial_i$, 
$B_i$ the image of $\partial_{i+1}$ and $H_i$ the homology group $Z_i / B_i$.
We say that $C_*$ is {\it acyclic\/} if $H_i$ vanishes for every $i$.

Let $c^i$ be a basis of $C_i$ and $c$ be the collection $\{c^i\}_{i \geq 0}$.
We call the pair $(C_*, c)$ a {\it based chain complex\/}, 
$c$ the preferred basis 
of $C_*$ and $c^i$ the preferred basis of $C_i$.
Let $h^i$ be a basis of $H_i$.

We construct another basis as follows. 
By the definitions of $Z_i$, $B_i$ and $H_i$, 
the following two split exact sequences exist.
\[
    0 \to Z_i \to C_i \stackrel{\partial_i}{\to} B_{i-1} \to 0,
\]
\[
    0 \to B_i \to Z_i \to H_i \to 0.
\]
Let ${\widetilde B}_{i-1}$ be a lift of $B_{i-1}$ to $C_i$ and 
${\widetilde H}_i$ a lift of $H_i$ to $Z_i$.
Then we can decompose $C_i$ as follows.
\begin{align*}
 C_i &= Z_i \oplus {\widetilde B}_{i-1} \\
     &= B_i \oplus {\widetilde H}_i  \oplus {\widetilde B}_{i-1} \\
     &= \partial_{i+1} {\widetilde B}_i \oplus {\widetilde H}_i  
                        \oplus {\widetilde B}_{i-1}.
\end{align*}
We choose $b^i$ a basis of $B_i$. 
We write $\tilde b^i$ for a lift of $b^i$ and 
 $\tilde h^i$ for a lift of $h^i$.
By the construction, 
the set $\partial_{i+1}(\tilde b^i) \cup \tilde h^i \cup \tilde b^{i-1}$ forms 
another ordered basis of $C_i$.
We denote simply this new basis 
by $\partial_{i+1}( \tilde b^i )\tilde h^i \tilde b^{i-1}$.
Then the definition of ${\rm tor}(C_*, c, h)$ is as follows.
\[
 {\rm tor}(C_*, c, h)
 =
 \prod_i^n [\partial_{i+1} ( \tilde b^i ) \tilde h^i \tilde b^{i-1} / c^i]^{(-1)^{i+1}} \in \F^*. 
\]
It is well known that ${\rm tor}(C_*, c, h)$ is independent of 
the choices of $\{b^i\}_{i \geq 0}$,
the lifts $\{\tilde b^i\}_{i \geq 0}$ and $\{\tilde h^i\}_{i \geq 0}$.

We also define the torsion ${\rm Tor}(C_*, c, h)$ with the sign term $(-1)^{|C_*|}$ 
as follows \cite{Turaev1}
\[
{\rm Tor}(C_*, c, h)
=
(-1)^{|C_*|}\cdot
{\rm tor}(C_*, c, h).
\]
Here 
\[
|C_*| = \sum_{i \geq 0} \alpha_i(C_*) \cdot \beta_i(C_*), 
\]
where
$\alpha_i(C_*) = \sum_{k=0}^{i} \dim C_k$
and 
$\beta_i(C_*) = \sum_{k=0}^i \dim H_k$.

\subsection{
Twisted chain complex and twisted cochain complex for CW-complex
}
\label{twisted_complex}
Let $W$ be a finite connected CW-complex and 
$\widetilde W$ its universal covering with 
the induced CW-structure.
Since the fundamental group $\pi_1(W)$ acts 
on $\widetilde W$ by the covering transformation,
the chain complex $C_* ({\widetilde W}; \Z)$ has a natural structure 
of a left $\Z[\pi_1(W)]$-module.
We denote by $\rho$ a homomorphism from $\pi_1(W)$ to $G$.
We regard the Lie group $\mathfrak{g}$ as a right $\Z[\pi_1(W)]$-module by
$\mathfrakg \times \pi_1(W) \ni (v, \gamma) \mapsto 
 Ad_{\rho(\gamma^{-1})}(v) \in \mathfrakg$.
We use the notation $\grho$ for $\mathfrakg$
with the right $\Z[\pi_1(W)]$-module structure.
Following \cite{KL, Porti}, we introduce the following notations.
Set
\begin{align*}
C_*(W ; \grho) &= \mathfrak{g} \otimes_{Ad \circ \rho} C_*(\widetilde W ; \Z),\\
C_*(W; \gtilde_{\rho}) &=
\mathfrakg(t) \otimes_{\alpha \otimes Ad\circ \rho} C_*(\widetilde W; \Z)
\end{align*}
where $\mathfrakg (t)$ is $\F(t) \otimes \mathfrakg$ and 
$\alpha$ is a surjective homomorphism from $\pi_1(W)$ to the multiplicative
group $\langle t \rangle$.
Note that 
$f \otimes v \otimes (\gamma\cdot \sigma) 
= f \cdot t^{\alpha(\gamma)} \otimes Ad_{\rho(\gamma^{-1})}(v)\otimes \sigma$.
We call $C_*(W; \grho)$ {\it the $\grho$-twisted chain complex\/}
and $C_*(W; \gtilde_{\rho})$ {\it the $\gtilde_{\rho}$-twisted chain complex\/} 
of $W$.
We also denote by $C^*(W; \grho)$ 
the $\F$-module consisting of the $\pi_1(W)$-equivalent homomorphisms 
from $C_*(\widetilde W; \Z)$ to $\mathfrakg$, 
i.e., a homomorphism $h$ satisfies 
$h(\gamma\cdot \sigma) = h(\sigma)\cdot \gamma^{-1}$ for 
$\gamma \in \pi_1(W)$.
We call $C^*(W; \grho)$ {\it the $\grho$-twisted cochain complex\/} of $W$.
$H_*(W ; \grho)$ and $H^* (W; \grho)$ denote 
the homology and cohomology groups of the $\grho$-twisted chain and cochain 
complexes.  

\subsection{The Reidemeister torsion for twisted chain complex}
We keep the notation of the previous subsection.
Let $e^{(i)}_1, \ldots, e^{(i)}_{n_i}$ be the set of $i$-dimensional cells of $W$.
We take a lift $\tilde e^{(i)}_{j}$ of the cell $e^{(i)}_j$ in $\widetilde W$.
Then, for each $i$, 
$\tilde c^i =\{\tilde e^{(i)}_1, \ldots, \tilde e^{(i)}_{n_i}\}$ is a basis of 
the $\Z[\pi_1(W)]$-module $C_i({\widetilde W};\Z)$.
Let $\mathbf B= \{\mathbf a, \mathbf b, \mathbf c\}$ be a basis of $\mathfrakg$.
Then we obtain the following basis of $C_i(W;\grho)$:
\[
\mathbf c_{\mathbf B} 
=
 \{ \ldots , 
     \mathbf a \otimes \tilde e^{(i)}_1, 
     \mathbf b \otimes \tilde e^{(i)}_1, 
     \mathbf c \otimes \tilde e^{(i)}_1,
    \ldots ,
     \mathbf a \otimes \tilde e^{(i)}_{n_i},
     \mathbf b \otimes \tilde e^{(i)}_{n_i},
     \mathbf c \otimes \tilde e^{(i)}_{n_i},
     \ldots 
  \}.
\]
When $\mathbf h^i = \{h^i_1, \ldots, h^i_{k_i}\}$ is 
a basis of $H_i(W;\grho)$, 
we denote by $\mathbf h$ the basis 
$\{\mathbf h^0, \ldots, \mathbf h^{\dim W}\}$ of $H_*(W; \grho)$.
Then ${\rm Tor}(C_*(W;\grho), \mathbf c_{\mathbf B}, \mathbf h) \in \F^*$ 
is well defined.
Furthermore adding a sign-refinement term into 
${\rm Tor}(C_*(W;\grho), \mathbf c_{\mathbf B}, \mathbf h)$, 
we define {\it the Reidemeister torsion\/} of $(W, \rho)$
 as a vector in some $1$-dimensional vector space
as follows.

\begin{definition}[\cite{dubois1, dubois2}]
Let $c_{\R}$ be the basis over $\R$ of $C_*(W; \R)$.
Choose an orientation $\mathfrak o$ of the real vector space 
$\oplus_{i\geq 0}H_i(W;\R)$
and provide $H_*(W; \R)$ with a basis 
$h_{\mathfrak o}=\{h^0, \ldots, h^{\dim W}\}$ 
such that each $h^i$ is a basis of $H_i(W; \R)$ and 
the orientation determined by $h_{\mathfrak o}$ agrees with $\mathfrak o$.
Let $\tau_0$ be either $+1$ or $-1$ according to the sign of 
${\rm Tor}(C_*(W;\R), c_{\R}, h_{\mathfrak o})$. 
Then
we define the Reidemeister torsion 
$\T(W, \grho, \mathfrak o)$ by 
\[
\T(W, \grho, \mathfrak o)
=
\tau_0 \cdot {\rm Tor}(C_*(W;\grho), \mathbf c_{\mathbf B}, \mathbf h)
\otimes_{i \geq 0} \det \mathbf h^i 
\in  Det\, H_*(W;\grho),
\]
where $\det \mathbf h^i = h^{(i)}_1 \wedge \ldots \wedge h^{(i)}_{k_i}$ and 
$
Det\, H_*(W;\grho)
=
\otimes_{i=0}^{\dim W} (\wedge^{\dim H_i} H_i(W; \grho))^{(-1)^{i}}.
$
Here $V^{-1}$ means the dual space of a vector space $V$ and 
the dual basis of 
$\det \mathbf h^i = h^{(i)}_1 \wedge \ldots \wedge h^{(i)}_{k_i}$
is $ h^{(i) *}_1 \wedge \ldots \wedge h^{(i) *}_{k_i}$
where $ h^{(i) *}_j$ is the dual element of $h^{(i)}_j$.
\end{definition}

We made some choices in the definition of $\T(W, \grho, \mathfrak o)$.
However the following well-definedness is known \cite[p.10]{Porti}:
\begin{itemize}
\item
The sign of $\T(W, \grho, \mathfrak o)$ is determined by 
the homology orientation $\mathfrak o$
i.e., if we choose the other homology orientation, 
then the sign of $\T(W, \grho, \mathfrak o)$ changes;
\item
$\T(W, \grho, \mathfrak o)$ does not depend on the choice of 
the lift $\tilde e^{(i)}_j$ for each cell $e^{(i)}_j$; 
\item
$\T(W, \grho, \mathfrak o)$ does not depend on 
the choice of the basis $\mathbf h$ in $\oplus_{i\geq 0} H_i(W;\grho)$.
\end{itemize}

We also have the following well-definedness.
\begin{lemma}
If the Euler characteristic of $W$ is equal to zero, 
then $\T(W, \grho, \mathfrak o)$ does not depend on 
the choice of the basis of $\mathfrakg$.
\end{lemma}

\begin{proof}
This follows from the definition.
\end{proof}

Similarly we define the Reidemeister torsion of 
the twisted $\gtilde_{\rho}$-chain complex.

\begin{definition}
We define $\T(W, \gtilde_{\rho}, \mathfrak o)$
by
\[
\T(W, \gtilde_{\rho}, \mathfrak o)
=
\tau_0 \cdot
{\rm Tor}(C_*(W; \gtilde), \mathbf  1 \otimes c_{\mathbf B}, \mathbf h)
\otimes_{i\geq 0} \det \mathbf h^i.
\]
\end{definition}
$\T(W, \gtilde_{\rho}, \mathfrak o)$ has the indeterminacy of $t^m$ 
where $m \in \Z$.
This indeterminacy is caused by the choice of the lifts $\{\tilde e^{(i)}_j\}$
and the action of $\alpha$.

\medskip

It is also known that the sign refined torsion
$\tau_0 \cdot {\rm Tor}(C_*(W; \grho), \mathbf c_{\mathbf B}, \mathbf h)$
has the invariance under simple homotopy equivalences,
and that it satisfies the following {\it Multiplicativity property\/}.
%
%
Suppose we have the following exact sequence of based chain complexes:
\begin{equation}\label{m_prop}
0 \to (C'_*, c') \to (C_*, c'\cup \bar c'') \to (C''_*, c'') \to 0
\end{equation}
where these chain complexes are 
based chain complexes which consist of vector spaces with bases.
Here we denote 
bases of $C'_*, C''_*$ by $c', c''$ and a lift of $c''$ to $C_*$ 
by $\bar c''$.
For each $i$, 
fix the volume forms on  $C'_i, C_i, C''_i$ by using given bases
and choose volume forms on $H_i(C'_*), H_i(C_*)$ and $H_i(C''_*)$. 
There exists the long exact sequence in homology 
associated to the short exact sequence $(\ref{m_prop})$:
\[
 \cdots \to H_i(C'_*) \to H_i(C_*) \to H_i(C''_*) \to H_{i-1}(C'_*) \to \cdots .
\]
We denote by $\mathcal H_*$ this acyclic complex.
Note that this acyclic complex is a based chain complex.

\begin{proposition}[Multiplicativity property \cite{Milnor66, Turaev2}]
\label{M_property}
We have
\[
{\rm Tor}(C_*) 
= 
(-1)^{\alpha(C'_*, C''_*)+\varepsilon(C'_*, C_*, C''_*)}
{\rm Tor}(C'_*)\cdot {\rm Tor}(C''_*)\cdot {\rm tor}(\mathcal H_*),
\]
where
\begin{align*}
\alpha(C'_*, C''_*)
&=
\sum_{i\geq 0} \alpha_{i-1}(C'_*) \alpha_i(C''_*) \in \Z/ 2\Z, \\
\varepsilon(C'_*, C_*, C''_*)
&=
\sum_{i\geq 0}
[(\beta_i(C_*) +1)(\beta_i(C'_*)+\beta(C''_*)) + \beta_{i-1}(C'_*)\beta(C''_*)]
\in \Z/2\Z.
\end{align*}
\end{proposition}

\subsection{On the representation spaces.}
Let $\pi$ be a finitely generated group and
we denote by $R(\pi, G)$ the space of $G$-representations of $\pi$.
We define the topology of this space by compact-open topology.
Here we assume that $\pi$ has the discrete topology and 
the Lie group $G$ has the usual one.
A representation $\rho:\pi \to G$ is called {\it central\/}
if $\rho(\pi) \subset \{\pm \mathbf 1\}$.

A representation $\rho$ is called {\it abelian\/}
if its image $\rho(\pi)$ is an abelian subgroup of $G$.
A representation $\rho$ is called {\it reducible\/}
if there exists a proper non-trivial subspace $U$ of $\C^2$ such that
$\rho(g)(U) \subset U$ for any $g \in \pi$.
A representation $\rho$ is called {\it irreducible\/} if it is not reducible.
We denote by $R^{red}(\pi, G)$ the subset of reducible representations and
by $R^{irr}(\pi, G)$ the subset of irreducible ones.
Note that all abelian representations are reducible.
The Lie group $G$ acts on $R(\pi, G)$ by conjugation.
We write $[\rho]$ for the conjugacy class of $\rho \in R(\pi, G)$,
and we denote by $\hat R(\pi, G)$ the quotient space $R(\pi, G)/ G$.  

\medskip

If $G$ is $\SU(2)$, then
one can see that the reducible representations are exactly abelian ones.
Note that this does not hold for the case of $\SL(2, \C)$-representations.
The action by conjugation of $\SU(2)$ on $R(\pi, \SU(2))$ factors through 
$\SO(3) = \SU(2)/\{\pm \mathbf 1\}$.
This action is free on the $R^{irr}(\pi, \SU(2))$.
We set $\hat R^{irr}(\pi, \SU(2)) = R^{irr}(\pi, \SU(2))/\SO(3)$.

If $G$ is $\SL(2, \C)$, then
the quotient space $\hat R(\pi, \SL(2, C))$ is not Hausdorff in general.
Following \cite{MS},
we will focus on the {\it character variety\/} 
$X(\pi; \SL(2, \C))$ 
which is the set of ${\it characters\/}$ of $\pi$.
Associated to the representation $\rho \in R(\pi, \SL(2, \C))$,
its character $\chi_{\rho}: \pi \to \C$, defined by 
$\chi_{\rho}(g)=\trace(\rho(g))$.
In some sense, $X(\pi, \SL(2, \C))$ is 
the ``algebro quotient'' of $R(\pi, \SL(2, \C))$ 
by ${\rm PSL}(2, \C)$.
It is well known that 
$R(\pi, \SL(2, \C))$ and $X(\pi)$ have the structure of complex algebraic 
affine sets 
and 
two irreducible representations of $\pi$ in $\SL(2, \C)$ with 
the same character are conjugate by an element of $\SL(2, \C)$.
(For the details, see \cite{MS}.)

\subsection{The Reidemeister torsion for knot exteriors.}
\label{Review_torsion_knot_ex}

In this subsection, we recall $\lambda$-regular representations and 
how to construct distinguished bases of $\grho$-twisted homology groups 
of knot exteriors 
for a $\lambda$-regular representation $\rho$.
These definitions have originally been given in \cite{Porti}.
The original definitions are written 
in terms of the $\grho$-twisted cohomology group.
We introduce the homology version by using the duality between
the twisted homology and cohomology associated to {\it the Kronecker pairing\/}
$C_*(W;\grho) \times C^*(W;\grho) \ni (\xi \otimes \sigma, v)
\mapsto (v(\sigma), \xi)_{\mathfrakg} \in \F$
\cite[p. 11]{Porti}.

\medskip

Let $K$ be a knot in a homology three sphere $M$.
We give a knot exterior $M_K$
the canonical homology orientation defined as follows.
It is well known that the $\R$-vector space
\[
H_*(M_K; \R) = H_0(M_K; \R) \oplus H_1(M_K; \R)
\]
has the basis $\{ [pt], [\mu]\}$.
Here $[pt]$ is the homology class of a point and 
$[\mu]$ is the homology class of a meridian of $K$.
We denote by $\mathfrak o$ the orientation induced by $\{[pt], [\mu]\}$.

We calculate the twisted homology groups of a circle and 
a $2$-dimensional torus 
before giving the definition of a natural basis of $H_*(M_K; \grho)$.
Here $S^1$ consists of one $0$-cell $e^{(0)}$ and one $1$-cell $e^{(1)}$.
\begin{lemma}\label{Lemma:homology_circle}
Suppose that $G$ is $\SU(2)$.
If $\rho \in R(\pi_1(S^1), G)$ is central, 
then $H_*(S^1; \grho)=\mathfrakg \otimes H_*(S^1; \R)$.
If $\rho$  is non-central, then we have
\begin{align*}
H_1(S^1; \grho) &= \R [P_{\rho} \otimes \tilde e^{(1)}], \,and \\
H_0(S^1; \grho) &= \R [P_{\rho}  \otimes \tilde e^{(0)}]
\end{align*}
where $P_{\rho}$ is a vector in $\mathfrakg$, 
which satisfies that $Ad(\rho(\gamma))(P_{\rho})=P_{\rho}$
for any $\gamma \in \pi_1(S^1)$.

Suppose that $G$ is $\SL(2, \C)$.
If $\rho \in R(\pi_1(S^1), G)$ is central, 
then $H_*(S^1; \grho)=\mathfrakg \otimes H_*(S^1; \C)$.
If $\rho$  is non-central and $\rho(\pi_1(S^1))$ has no parabolic elements, 
then we have
 \begin{align*}
 H_1(S^1; \grho)&=\C [P_{\rho} \otimes \tilde e^{(1)}], \,and\\
 H_0(S^1; \grho)&=\C [P_{\rho} \otimes \tilde e^{(0)}]
 \end{align*}
 where $P_{\rho}$ is a vector in $\mathfrakg$, 
which satisfies that $Ad(\rho(\gamma))(P_{\rho})=P_{\rho}$
for any $\gamma \in \pi_1(S^1)$.
If $\rho$  is non-central and the subgroup $\rho(\pi_1(S^1))$ is contained 
in a subgroup which consists of parabolic elements, 
then we have
\[
 H_1(S^1; \grho)=\C [P_{\rho} \otimes \tilde e^{(1)}].
\]
\end{lemma}
\begin{proof}
This is a consequence of the following fact of homology of groups.
For $G = \Z$, it follows that
$
  H_0(G;N) = H^1(G;N) = N_{G} 
$
and 
$
  H^0(G;N) = H_1(G;N) = N^{G}
$
where $G$ is a group, $N$ is a $N$-module,
$N_{G}$ is the group of invariants of $N$ and
$N^{G}$ is the group of co-invariants of $N$ 
(for the details, see \cite{Brown}).
\end{proof}

We denote by $T^2$ a $2$-dimensional torus.
Here $T^2$ consists of 
one $0$-cell $e^{(0)}$, two $1$-cells $e^{(1)}_1, e^{(1)}_2$ and 
one $2$-cell $e^{(2)}$.
We denote each cell $e^{(0)}, e^{(1)}_1, e^{(1)}_2$ and $e^{(2)}$ 
by $pt, \mu, \lambda$ and $T^2$. 
One can also calculate 
the $\grho$-twisted homology groups of $C_*(T^2; \grho)$ as follows.

\begin{lemma}\label{Lem:homology_torus}
Suppose that $G$ is $\SU(2)$. 
If $\rho \in R(\pi_1(T^2), G)$ is central, 
then $H_*(T^2; \grho)=\mathfrakg \otimes H_*(T^2; \R)$.
If $\rho \in R(\pi_1(T^2), G)$ is non-central, then we have
\begin{align*}
H_2(T^2; \grho)
&=\R [\, P_{\rho} \otimes \widetilde T^2], \\
H_1(T^2; \grho)
&=\R [\, P_{\rho} \otimes \tilde \mu\,] \oplus 
  \R [\, P_{\rho} \otimes \tilde \lambda\,],\\
H_0(T^2; \grho) &=\R [\,P_{\rho} \otimes \widetilde pt\,]
\end{align*}
where $P_{\rho}$ is a vector of $\mathfrakg$ 
such that $Ad_{\rho(\gamma)}(P_{\rho}) = P_{\rho}$ 
for any $\gamma \in \pi_1(T^2)$.

Suppose that $G$ is $\SL(2, \C)$.
If $\rho \in R(\pi_1(T^2), G)$ is central, 
then $H_*(T^2; \grho)=\mathfrakg \otimes H_*(T^2; \C)$.
If $\rho \in R(\pi_1(T^2), G)$ is non-central and 
 $\rho(\pi_1(T^2))$ contains a non-parabolic element, 
then we have
\begin{align*}
H_2(T^2; \grho)
&=\C [\, P_{\rho} \otimes \widetilde T^2],\\
H_1(T^2; \grho)
&=\C [\, P_{\rho} \otimes \tilde \mu\,] \oplus 
  \C [\, P_{\rho} \otimes \tilde \lambda \,],\\
H_0(T^2; \grho) &=\C [\, P_{\rho} \otimes \widetilde pt \,] 
\end{align*}
where $P_{\rho}$ is a vector of $\mathfrakg$ 
such that $Ad_{\rho(\gamma)}(P_{\rho}) = P_{\rho}$ 
for any $\gamma \in \pi_1(T^2)$.

If $\rho \in R(\pi_1(T^2), G)$ is non-central and 
the subgroup $\rho(\pi_1(T^2))$ is contained in a subgroup 
which consists of parabolic elements,
then we have
\[
H_2(T^2; \grho)=\C [P_{\rho} \otimes \widetilde T^2] 
\]
and $[P_{\rho} \otimes \tilde \lambda]$ is a non-zero class in $H_1(M_K; \grho)$.
\end{lemma}

\begin{proof}
This is a consequence of \cite[Proposition 3.18]{Porti}.
\end{proof}

Next we give the definition of regular representations for $\pi_1(M_K)$
in terms of the twisted $\grho$-chain complex.

\begin{definition}[regular representations {\cite[p.83]{Porti}}] 
We say that $\rho$ is {\it regular\/}
if $\rho$ is irreducible and $\dim_{\F} H_1(M_K; \grho)=1$.

We let $\gamma$ be a simple closed curve in $\partial M_K$.
We say that $\rho$ is {\it $\gamma$-regular\/}
if :
\begin{itemize}
\item[(1)] $\rho$ is regular;
\item[(2)] an inclusion $\iota:\gamma \hookrightarrow M_K$ 
           induces the surjective homomorphism
\[
\iota_* : H_1(\gamma; \grho) \to H_1(M_K; \grho); \quad  and
\]
 \item[(3)] if $\trace(\rho(\pi_1(\partial M_K))) \subset \{\pm 2\}$, 
  then $\rho(\gamma) \not = \pm \mathbf 1$.
\end{itemize}
\end{definition}

We fix an invariant vector $P_{\rho} \in \mathfrakg$ as above.
Let $\gamma$ be a simple closed curve in $\partial M_K$.
An inclusion $\iota: \gamma \hookrightarrow M_K$ and 
the the Kronecker pairing between homology and cohomology induce 
the linear form $f^{\rho}_{\gamma} : H^1(M_K; \grho) \to \F$.
By Lemma \ref{Lemma:homology_circle}, it is explicitly described by  
\[
f^{\rho}_{\gamma}(v)
=
( \iota_*([\tilde \gamma \otimes P_{\rho} ]), v )
=
(P_{\rho}, v(\tilde \gamma))_{\mathfrakg}
 \quad 
\text{for any}\,\, 
v \in H^1(M_K; \grho).
\] 
An alternative formulation of 
$\gamma$-regular representations is given in \cite{dubois2, Porti}.
Similarly, we can also give the following alternative formulation of 
the $\gamma$-regularity in our conventions.

\begin{proposition}\label{check_gamma_regular}
A representation $\rho \in R^{irr}(\pi_1(M_K), G)$ is $\gamma$-regular
if and only if 
the linear form $f^{\rho}_{\gamma}: H^1(M_K; \grho) \to \F$ is an isomorphism.
\end{proposition}

\begin{proof}
If $f^{\rho}_{\gamma}$ is an isomorphism,
then we have that $\dim_{\F}H^1(M_K; \grho)=1$ and 
$\iota_*([P_{\rho} \otimes \tilde \gamma])$ is a non-zero class 
in $H_1(M_K; \grho)$.
It follows from 
the Kronecker pairing between the $\grho$-twisted homology and cohomology 
that $\dim_{\F}H_1(M_K; \grho)$ is also one.
Hence $\iota_*$ is surjective.
If $\rho$ is $\gamma$-regular, then
we have that $\dim_{\F}H_1(M_K; \grho)=1$ and 
$\iota_* : H_1(\gamma; \grho) \to H_1(M_K; \grho)$ is surjective. 
We denote a generator of $H_1(M_K; \grho)$ by $\sigma$.
There exists an element $[v \otimes \tilde \gamma]$ of $H_1(\gamma; \grho)$ 
such that $\iota_*([v \otimes \tilde \gamma])=\sigma$.

If $\rho(\gamma)$ is central, then $v$ satisfies that 
$Ad(\rho(\gamma'))(v) = v$ for any $\gamma' \in \pi_1(\partial M_K)$.
Therefore $\iota_*([v \otimes \tilde \gamma ])$ induces 
the isomorphism $f^{\rho}_{\gamma}$.

Suppose that $\rho(\gamma)$ is non-central, 
then $H_1(\gamma; \grho)$ is generated by $[P_{\rho} \otimes \tilde \gamma\,]$.
There exists an element $c \in \F^*$ 
such that $[v \otimes \tilde \gamma] = c [P_{\rho} \otimes \tilde \gamma]$.
Hence $\iota_*([P_{\rho} \otimes \tilde \gamma])$ is a non-zero class 
in $H_1(M_K; \grho)$.
Therefore $\iota_*([P_{\rho} \otimes \tilde \gamma])$ induces 
the isomorphism $f^{\rho}_{\gamma}$.
\end{proof}

We define a reference generator of $H_1(M_K; \grho)$ 
by using the above isomorphism $f^{\rho}_{\gamma}$.

Let $\rho$ be a $\lambda$-regular representation of $\pi_1(M_K)$.
By Lemma \ref{Lem:homology_torus}, 
the reference generator of $H_1(M_K; \grho)$ is defined by 
\[
h^{(1)}_{\rho}(\lambda) = \iota_*([P_{\rho} \otimes \tilde \lambda\, ]).
\]

\medskip 

Moreover the reference generator of $H_2(M_K; \grho)$ is defined as follows.

\begin{lemma}[Cor.\,3.23 \cite{Porti}]
Let $i : \partial M_K \hookrightarrow M_K$ be an inclusion map.
If $\rho \in R(\pi_1(M_K), G)$ is $\gamma$-regular, 
then
we have the isomorphism $i_* : H_2(\partial M_K; \grho) \to H_2(M_K; \grho)$.
\end{lemma}

Using this isomorphism $i_*$, 
we define the reference generator of $H_2(M_K; \grho)$ by
\[
  h^{(2)}_{\rho} = i_* ([P_{\rho} \otimes \widetilde{\partial M_K}]).
\]

\begin{remark}
The reference generators of $H^1(M_K; \grho)$ and $H^2(M_K; \grho)$ 
have been defined in \cite{dubois1, dubois2, Porti} by using another 
metric of $\mathfrakg$.
If we define reference generators of $H^1(M_K; \grho)$ and $H^2(M_K; \grho)$ 
by using our metric $(\, , \,)_{\mathfrakg}$, 
then the resulting generators become 
the dual bases of $h^{(1)}_{\rho}(\lambda)$ and $h^{(2)}_{\rho}$
from the above propositions.  (For the details, see \cite{dubois2, Porti}.)
\end{remark}

\noindent
We recall the definition of {\it the twisted Reidemeister torsion\/} for knot exteriors.
Let $\rho:\pi_1(M_K) \to G$ be a $\lambda$-regular representation.
We define ${\mathbb T}^{K}_{\rho}$ by 
the coefficient of the Reidemeister torsion $\T(M_K, \grho, \mathfrak o)$
where we choose the reference generators $h^{(1)}_{\rho}(\lambda), h^{(2)}_{\rho}$ 
as a basis of $H_*(M_K; \gtilde)$, i.e.,
$\mathbb T^K_{\lambda}$ is given explicitly by
\[
{\mathbb T}^{K}_{\lambda}(\rho)
=
\tau_{0}\cdot
{\rm Tor}(C_{*}(M_K; \grho), \mathbf c_{\mathbf B}, 
    \{h^{(1)}_{\rho}(\lambda), h^{(2)}_{\rho}\})
\in \F^{*}.
\]

Given the reference generator of $H_*(M_K; \grho)$, 
the basis of the determinant line $Det\, H_*(M_K; \grho)$ is also given.
This means that a trivialization of the line bundle $Det\, H_*(M_K; \grho)$ 
at $\rho$ is given. 
The Reidemeister torsion $\T(M_K, \grho, \mathfrak o)$ is 
a section of the line bundle $Det\, H_*(M_K; \grho)$.
We can regard ${\mathbb T}^{K}_{\lambda}$ 
as a section of the line bundle $Det\, H_*(M_K; \grho)$ 
over $\lambda$-regular representations 
with respect to the trivialization 
by $\{h^{(1)}_{\rho}(\lambda), h^{(2)}_{\rho}\}$.
We also call $\mathbb T^K_{\lambda}$ {\it the twisted Reidemeister torsion\/}.

\section{A relationship  between acyclic Reidemeister torsion and
               non-acyclic Reidemeister torsion}
\label{Main_theorem}
\subsection{The statement of main theorem}
Our purpose is to express the twisted Reidemeister torsion 
by using a limit of the acyclic Reidemeister torsion.

Let $K$ be a knot in a homology three sphere $M$ and $M_K$ its exterior.
One of the invariants which we will investigate is 
the twisted Reidemeister torsion ${\mathbb T}^K_{\lambda}$.
The other is 
the acyclic Reidemeister torsion $\T(M_K, \gtilde_{\rho}, \mathfrak o)$.
This invariant coincides with the twisted Alexander invariant of $\pi_1(M_K)$
\cite{Kitano}.
The twisted Alexander invariant is computed by using 
the Fox calculus \cite{KL, Kitano}.
We prove that the twisted Reidemeister torsion may be expressed as 
the differential coefficient of the twisted Alexander invariant of $\pi_1(M_K)$.

The invariant $\T(M_K, \gtilde_{\rho}, \mathfrak o)$ is only defined 
when the local system $C_*(M_K; \gtilde_{\rho})$ is acyclic.
On the other hand, the twisted Reidemeister torsion ${\mathbb T}^K_{\lambda}$
is defined on the set of $\lambda$-regular representations of $\pi_1(M_K)$.
We need to check 
whether the local system $C_*(M_K; \gtilde_{\rho})$ is acyclic
for a $\lambda$-regular representation $\rho$.

\begin{proposition}\label{PropositionA}
Let $\rho$ be an $\SU(2)$ or $\SL(2, \C)$-representation of a knot group.
If $\rho$ is $\lambda$-regular, 
then the twisted chain complex $C_*(M_K;\widetilde{\mathfrak{g}}_\rho)$ 
is acyclic. 
\end{proposition}

Note that for a knot exterior in a homology $3$-sphere, 
the homomorphism $\alpha$ satisfies $\alpha(\mu)=t$ 
where $\mu$ is the meridian of the knot.

\medskip

Therefore 
${\mathbb T}^K_{\lambda}$ and $\T(M_K, \gtilde_{\rho}, \mathfrak o)$
are well defined on $\lambda$-regular representations.
By the definitions, 
the twisted Reidemeister torsion ${\mathbb T}^K_{\lambda}$ is 
an element of $\F^*$ and 
the twisted Alexander invariant $\T(M_K, \gtilde_{\rho}, \mathfrak o)$
is an element of $\F(t)^*$.
Actually
the following relation 
between ${\mathbb T}^K_{\lambda} \in \F^*$ and the rational function 
$\T(M_K, \gtilde_{\rho}, \mathfrak o) \in \F(t)^*$.

\begin{theorem}\label{thm:main_theorem}
If $\rho$ is a $\lambda$-regular representation, then
the acyclic Reidemeister torsion $\T(M_K, \gtilde_{\rho}, {\mathfrak o})$ 
for $\rho$ 
has a simple zero at $t=1$. 
Moreover the following holds:
\[
  {\mathbb T}_{\lambda}^{K}(\rho)
  = - \lim_{t \to 1}\frac{\T(M_K, \gtilde_{\rho}, {\mathfrak o})(t)}{t-1}
  = - \left. \frac{d}{dt} \T(M_K, \gtilde_{\rho}, {\mathfrak o}) \right|_{t=1}. 
\]
\end{theorem}

This says that
we can compute the twisted Reidemeister torsion ${\mathbb T}^K_{\lambda}$ 
algebraically
by using Fox calculus of the twisted Alexander invariant of $K$.

\subsection{Proof of Proposition \ref{PropositionA}}

We prove Proposition \ref{PropositionA} by using 
the $\lambda$-regularity of $\rho$.

\begin{proof}[Proof of Proposition \ref{PropositionA}]
It is well known that any compact connected triangulated $3$-manifold 
whose boundary is non-empty 
and consists of tori can be collapsed into a $2$-dimensional sub-complex
(see II. Cor. 11.9 in \cite{Turaev1}).
Moreover, by the simple-homotopy extension theorem, 
every CW-complex has the simple-homotopy type of a CW-complex which has only one vertex.
We denote this $2$-dimensional CW-complex by $W$ 
and this deformation from $M_K$ to $W$ by $\varphi$.
Since two $\gtilde_{\rho}$-twisted homology groups $H_*(M_K; \gtilde_{\rho})$ 
and $H_*(W;\gtilde_{\rho})$ are isomorphic, 
we prove that $H_*(W; \gtilde_{\rho})$ vanishes in the following.

The fact that $H_0(W; \gtilde_{\rho}) = 0$ is proved in 
\cite[Proposition 3.5]{KL}.
Since the Euler characteristic of $W$ is zero, 
the dimension of $H_1(W; \gtilde_{\rho})$ is equal to 
that of $H_2(W; \gtilde_{\rho})$.
We must prove that the dimension of $H_2(W;\gtilde_{\rho})$ over $\F(t)$ is zero.
It is enough to prove that the rank over $\F[t, t^{-1}]$ 
of the second homology group of the following local system is zero:
\[
  C_*(W;\grho[t, t^{-1}]) 
  = 
  \mathfrakg[t, t^{-1}] \otimes_{\alpha \otimes Ad \circ \rho} C_*(\tilde W;\Z) 
\]
where $\mathfrakg[t, t^{-1}]$ is $\F[t, t^{-1}] \otimes \mathfrakg$.
We denote the homology group of this chain complex by 
$H_*(W;\grho[t, t^{-1}])$.
Suppose that the rank of $H_2(W;\grho[t, t^{-1}]) > 0$.

There exists the long exact homology sequence \cite{Spanier}:
\[
0 \to H_2(W;\grho[t, t^{-1}]) 
  \xrightarrow{(t-1)\cdot} H_2(W;\grho[t, t^{-1}])
  \xrightarrow{t=1} H_2(W;\grho)
  \xrightarrow{\Delta} H_1(W;\grho[t, t^{-1}])
  \to \cdots
\]
associated to the short exact sequence:
\[
 0 \to \mathfrakg[t, t^{-1}] 
   \xrightarrow{(t-1)\cdot} \mathfrakg[t, t^{-1}]
   \xrightarrow{t=1} \mathfrakg
   \to 0. 
\]
Since the rank of $H_2(W; \grho[t,t^{-1}])$ is not zero, 
the multiplication with $(t-1)$ is not surjective.
Hence the image of the evaluation map $(t=1)$ is not trivial and therefore
surjective since the dimension of $H_2(W; \grho)$ is only one.
This implies that $\Delta$ is trivial. 
On the other hand the equation
\[
\partial (1 \otimes P_{\rho} \otimes \widetilde{\varphi(\partial M_K)})
= (t-1)\cdot(1\otimes P_{\rho} \otimes \widetilde{\varphi(\lambda)})
\]
implies that 
$\Delta([P_{\rho}\otimes \widetilde{\varphi (\partial M_K)}]) 
= [1\otimes P_{\rho} \otimes \widetilde{\varphi(\lambda)}]$.
But $[1\otimes P_{\rho} \otimes \widetilde{\varphi(\lambda)}]$ can not be trivial
since it is mapped under the evaluation map $(t=1)$ to $[P_{\rho}\otimes \widetilde{\varphi(\lambda)}]$ and
the chain $P_{\rho} \otimes \widetilde{\varphi(\lambda)}$ represents
a non-zero homology class in $H_1(W; \grho)$.
This is a contradiction.
Therefore 
the rank of $H_2(W; \grho[t, t^{-1}])$  over $\F[t, t^{-1}]$ is zero.
Hence we have that $\dim_{\F(t)} H_2(W; \gtilde_{\rho})=0$. 
Also $\dim_{\F(t)} H_1(W; \gtilde_{\rho})$ is zero.
\end{proof}

\subsection{Proof of Theorem \ref{thm:main_theorem}}
At first, we prepare some notations and an algebraic proposition.

\medskip

Let $C_*$ is an $n$-dimensional chain complex 
which consists of left $G$-modules $M_i\, (1 \leq i \leq n)$ 
where $G$ is a group.
We denote by $C_*(V)$ the chain complex 
which consists of the vector spaces $V \otimes_{\rho} M_i$ 
where $V$ is a right $G$-vector space over $\F$ and 
$\rho$ is a homomorphism from $G$ to ${\rm Aut}(V)$.
Let $H_*(V)$ be the homology groups of $C_*(V)$,
$C'_*(V)$ the subchain complex 
which consists of a lift of $H_*(V)$ to $C_*(V)$
and $C''_*(V)$ the quotient of $C_*(V)$ by $C'_*(V)$.
We denote by $h(V), c'$ and $c''$
the bases of $H_*(V), C'_*(V)$ and $C''_*(V)$.
Note that $c'$ is a lift of $h(V)$ to $C_*(V)$.
If there exists a homomorphism $\alpha$
from $G$ to the multiplicative group $\langle t \rangle$,
we denote by $C_*(V(t))$ 
which consists of vector spaces $V(t) \otimes_{\alpha \otimes \rho} M_i$.
Here we denote $\F(t) \otimes V$ by $V(t)$.
Moreover let $C'_*(V(t))$ be the subchain complex which is given by 
extending the coefficients of $C'_*(V)$ to $\F(t)$ by using $\alpha$
and $C''_*(V(t))$ the quotient of $C_*(V(t))$ by $C'_*(V(t))$.
\begin{proposition}\label{prop:alg_preparation}
We assume that $C_*(V(t))$ and $C'_*(V(t))$ are acyclic.
The following relation holds:
\begin{eqnarray}
\lefteqn{\lim_{t \to 1}(-1)^{\alpha'}
             \frac{{\rm Tor}(C_*(V(t)),  1 \otimes c' \cup 1 \otimes \bar c'' )}
                  {{\rm Tor}(C'_*(V(t)), 1 \otimes c')}
}\\
&&= (-1)^{\varepsilon' + |C_*(V)|}{\rm Tor}(C_*(V), c' \cup \bar c'' , h(V))\nonumber
\end{eqnarray}
where $\bar c''$ is a lift of $c''$ to $C_*(V)$, 
$\alpha'$ is $\alpha(C'_*(V(t)), C''_*(V(t)))$
in Proposition \ref{M_property},
and $\varepsilon' \in \Z / 2\Z$ is given by
$\sum_{i=0}^{n-1} \dim_{\F} C''_i(V) \cdot \beta_i(C_*(V))$.
\end{proposition}

\begin{proof}
The chain complex $C''_*(V(t))$ is also acyclic 
from the long exact sequence of the pair $(C_*(V(t)), C'_*(V(t)))$. 
We can apply Proposition \ref{M_property} for the short exact sequence:
\[
  0 \to (C'_*(V(t)), 1 \otimes c') 
    \to (C_*(V(t)),  1 \otimes c' \cup 1 \otimes \bar c'') 
    \to (C''_*(V(t)), 1 \otimes c'')
    \to 0.
\]
Then, we obtain the following equation of the torsions.
\begin{equation}\label{eqn:result_M_property}
(-1)^{\alpha'} {\rm Tor}(C_*(V(t)), 1 \otimes c' \cup 1 \otimes \bar c'')
=
{\rm Tor}(C'_*(V(t)), 1 \otimes c') \cdot {\rm Tor}(C''_*(V(t)), 1 \otimes c'').
\end{equation}
Note that $\varepsilon(C'_*(V(t)), C_*(V(t)), C''_*(V(t)))=0$
because $C_*(V(t))$, $C'_*(V(t))$ and $C''_*(V(t))$ are acyclic.

Next we consider ${\rm Tor}(C''_*(V(t)), c'')$.
It follows 
from the long exact sequence of the pair $(C_*(V), C'_*(V))$ 
and the definition of $C'_*(V)$ 
that the chain complex $C''_*(V)$ is also acyclic.
Since $C''_*(V)$ is acyclic, 
we can choose a basis $\tilde b''^i$ of $\widetilde B''_i$ for each $i$. 
Here $\widetilde B''_i$ is a lift of 
$B''_{i} = {\rm Im}\, \partial_{i+1}(C''_{i+1}(V))$
to $C''_{i+1}(V)$.

\begin{claim}\label{claim:on_lift_b}
A subset $1 \otimes \tilde b''^i$ in $C''_{i+1}(V(t))$ 
generates a subspace on which the boundary operator $\partial_{i+1}$
is injective.
\end{claim}

\noindent
{\it Proof of Claim \ref{claim:on_lift_b}.\/}
If the determinant of the boundary operator 
restricted on $\F(t) \langle 1 \otimes \tilde b''^i \rangle$
is zero, 
then  substituting $1$ for the parameter $t$ 
we have that 
the determinant of the boundary operator 
restricted on $\F\langle \tilde b''^i \rangle$
is also zero.
This is a contradiction to the choices of $\tilde b''^i$.
\hfill (Claim \ref{claim:on_lift_b}) $\Box$

\medskip

Therefore ${\rm Tor}(C''_*(V(t)), 1 \otimes c'')$ is represented as 
\[
  \prod_{i=0}^{n}
  \left[
    \partial_{i+1} (1 \otimes \tilde b''^i) 1 \otimes \tilde b''^{i-1}
                    / 1 \otimes c''^{i}
  \right]^{(-1)^{i+1}}.
\]
We denote by $\tilde b^i$ 
a lift $1 \otimes \tilde b''^i$ to $C_*(V(t))$ simply.
Note that 
\begin{align*}
\lefteqn{
  \prod_{i=0}^n
  \left[
    \partial_{i+1}(1 \otimes \tilde b''^i)\, 1 \otimes \tilde b''^{i-1}
                   / 1 \otimes c''^i
  \right]^{(-1)^{i+1}}
}& \\
&=
  \prod_{i=0}^n
  \left[
    (1 \otimes c'^i) \, \partial_{i+1} (\tilde b^i) \, \tilde b^{i-1} 
    /  1 \otimes c'^i \cup 1 \otimes \bar c''^i
  \right]^{(-1)^{i+1}}.
\end{align*}

We substitute these results into 
the equation (\ref{eqn:result_M_property})
Then we have 
\begin{align}
\lefteqn{
\frac{{\rm Tor}(C_*(V(t)), 1 \otimes c' \cup 1 \otimes \bar c'')}
     {{\rm Tor}(C'_*(V(t)), 1 \otimes c')}
}
& \nonumber\\
&=
 {\rm Tor}(C''_*(V(t)), 1 \otimes c'') \nonumber \\
&= 
  \prod_{i=0}^n
  \left[
    (1 \otimes c'^i) \, \partial_{i+1} (\tilde b^i) \, \tilde b^{i-1} 
     / 1 \otimes c'^i \cup 1 \otimes \bar c''^i
  \right]^{(-1)^{i+1}} \nonumber\\
&=
  \prod_{i=0}^n
  (-1)^{\dim_{\F} B''_i \cdot \dim_{\F}H_i(V)}
  \left[
    \partial_{i+1} (\tilde b^i)\, (1 \otimes c'^i) \, \tilde b^{i-1} 
     / 1 \otimes c'^i \cup 1 \otimes \bar c''^i
  \right]^{(-1)^{i+1}}. \label{eqn:ratio_torsion}
\end{align}
The acyclicity of $C''_*(V)$ shows that
\[
\sum_{i=0}^n \dim_{\F} B''_i \cdot \dim_{\F}H_i(V)
\equiv 
\sum_{i=0}^{n-1} \dim_{\F}C''_i(V) \cdot \beta_i(C_*(V))
\quad ({\rm mod}\, 2).
\]
Substituting $1$ for $t$, 
the right hand side (\ref{eqn:ratio_torsion}) turns into 
\[
  (-1)^{\varepsilon'}
  \prod_{i=0}^n
  \left[
    \partial_{i+1} (\tilde b^i)\, \tilde h^i \, \tilde b^{i-1} 
     / c'^i \cup \bar c''^i
  \right]^{(-1)^{i+1}}.
\]
This is equal to $(-1)^{\epsilon' + |C_*(V)|}{\rm Tor}(C_*(V), c' \cup \bar c'' , h(V))$.
 
Although the left hand side is determined up to a factor $t^m (m \in \Z)$,
the limit at $t=1$ is determined 
because the factor $t^m$ does not affect taking a limit at $t=1$. 
\end{proof}
%
We can prove Theorem \ref{thm:main_theorem} 
as an application of Proposition \ref{prop:alg_preparation}.
\begin{proof}[Proof of Theorem \ref{thm:main_theorem}.]
As in the proof of Proposition \ref{PropositionA},
let $W$ be a $2$-dimensional CW-complex 
with a single vertex which has the same simple-homotopy type as $M_K$.
We denote the deformation from $M_K$ to $W$ by $\varphi$.
The compact $3$-manifold $M_K$ is simple homotopy equivalent to $W$.
It is enough to prove the theorem for $W$ because of the invariance of 
the simple homotopy equivalence for the Reidemeister torsion.
Let $\rho$ be a $\lambda$-regular representation of $\pi_1(M_K)$.
We denote by the same symbols $\rho$ and $\mathfrak o$ 
the representation  of $\pi_1(W)$
and the homology orientation of $H_*(W; \R)$
induced from that of $M_K$ under the map $\varphi$.

We define the subchain complex $C'_*(W; \grho)$ of 
the $\grho$-twisted chain complex $C_*(W; \grho)$ by
\[
C'_2(M_K; \grho) 
  = \F\langle P_{\rho} \otimes \widetilde{\varphi(\partial M_K)} \rangle,
\quad
C'_1(W; \grho) 
  = \F\langle  P_{\rho} \otimes\widetilde{\varphi(\lambda)} \rangle
\]
and $C_i(W; \grho) = 0\,(i \not = 1, 2)$
where $P_{\rho}$ is an invariant vector of $\mathfrakg$ 
such that
 $Ad_{\rho(\gamma)}(P_{\rho})=P_{\rho}$
for
any $\gamma \in \pi_1(\varphi(\partial M_K))$.
The modules of this subchain complex are lifts of homology groups 
$H_*(W; \grho)$.
By the definition, 
the boundary operators of $C'_*(W; \grho)$ are zero homomorphisms. 
Let $C''_*(W; \grho)$ be the quotient of $C_*(W; \grho)$ by $C'_*(W; \grho)$.
Similarly, 
we define the subcomplex $C'_*(W; \gtilde_{\rho})$ of $C_*(W; \gtilde_{\rho})$ 
to be 
\[
C'_2(W; \gtilde_{\rho}) 
 = \F(t)
   \langle 1 \otimes P_{\rho} \otimes \widetilde{\varphi(\partial M_K)}\rangle, 
\quad
C'_1(W; \gtilde_{\rho}) 
  = \F(t)
    \langle 1 \otimes P_{\rho} \otimes \widetilde{\varphi(\lambda)}\rangle
\]
and $C'_i(W) = 0$ for $i \not = 1, 2$.
The boundary operators of $C'_*(W; \gtilde_{\rho})$ is given by
\[
0  \to C'_2(W; \gtilde_{\rho}) 
   \xrightarrow{(t-1)\cdot} C'_1(W; \gtilde_{\rho}) 
   \to 0. 
\]
This shows that 
the subchain complex $C'_*(M_K; \gtilde_{\rho})$ is acyclic.
By Proposition \ref{PropositionA}, 
the $\gtilde_{\rho}$-twisted chain complex $C_*(M_K; \gtilde_{\rho})$ 
is also acyclic.

The twisted chain complex $C'_*(W; \grho)$ has the natural basis: 
\[
  c'
  =
  \{
    P_{\rho} \otimes \widetilde{\varphi(\partial M_K)},\,  
    P_{\rho} \otimes \widetilde{\varphi(\lambda)}
  \}.
\]
Let $c''$
be
a basis of $C''_*(W; \grho)$ 
and  
$\bar c''$
a lift of $c''$ to $C_*(W; \grho)$.
Applying Proposition \ref{prop:alg_preparation},
we have
\begin{eqnarray}\label{eqn:alg_prop_W}
\lefteqn{
\lim_{t \to 1}
  \frac{(-1)^{\alpha'} 
        {\rm Tor}(C_*(W; \gtilde_\rho), 1 \otimes c' \cup 1 \otimes \bar c'')}
       {{\rm Tor}(C'_*(W; \gtilde_\rho), 1 \otimes c')}
}\\
&&=
(-1)^{\varepsilon' + |C_*(W; \grho)|}{\rm Tor}(C_*(W; \grho), c' \cup \bar c'', \{h^{(1)}_\rho(\lambda), h^{(2)}_\rho\}).\nonumber
\end{eqnarray}
\begin{claim}\label{claim:proof_main_thm}\hphantom{ }\hfill
\begin{enumerate}
\item
${\rm Tor}(C'_*(W; \gtilde_\rho), 1 \otimes c' ) = t-1$. 
\item
$\alpha' \equiv 0\,({\rm mod}\,2)$.
\item
$\varepsilon' + |C_*(W; \grho)| \equiv 1 \, ({\rm mod}\, 2)$.
\end{enumerate}
\end{claim}
\noindent
{\it Proof of Claim \ref{claim:proof_main_thm}.\/}
$(1)$ It follows by the definition.
$(2)$ 
If we denote the number of $1$-cells of $W$ by $k$, 
the CW-complex $W$ has one $0$-cell, $k$ $1$-cells and $(k-1)$ $2$-cells.
We have 
$\alpha' 
= 0 \cdot (3k+2) + 1 \cdot (6k-2) + 2 \cdot (6k-2)
\equiv 0\, ({\rm mod}\, 2)$.
$(3)$
This follows from 
$\varepsilon' = (3k - 4) \cdot 1 \equiv 3k-4\, ({\rm mod}\, 2)$
and
$|C_*(W; \grho)|=
 3 \cdot 0 + (3k+3)\cdot 1 + (3k+3+3k-3)\cdot 2 
\equiv 3k+3 \, ({\rm mod} \, 2)
$.
\hfill (Claim \ref{claim:proof_main_thm})$\Box$

\medskip

The equation (\ref{eqn:alg_prop_W}) turns into
\[
  \lim_{t \to 1}
  \frac{{\rm Tor}(C_*(W; \gtilde_\rho), 1 \otimes c' \cup 1 \otimes \bar c'')}
       {t-1}
=
  -{\rm Tor}(C_*(W; \grho), c' \cup \bar c'', \{h^{(1)}_\rho(\lambda), h^{(2)}_\rho\}).
\]
Multiplying the both sides 
by the alternative products of the determinants of the base-change matrices 
\[
\prod_{i=0}^2
\left[
c'^i \cup \bar c''^i / \mathbf c_{\mathbf B}
\right]^{(-1)^{i+1}},
\]
we obtain the following equation:
\[
  \lim_{t \to 1}
  \frac{{\rm Tor}(C_*(W; \gtilde_\rho), \mathbf c_{\mathbf B})}
       {t-1}
=
  -{\rm Tor}(C_*(W; \grho),  \mathbf c_{\mathbf B}, \{h^{(1)}_\rho(\lambda), h^{(2)}_\rho\}).
\]
Finally multiplying the both sides by the sign $\tau_0$ gives
\[
  \lim_{t \to 1}
  \frac{\T(W, \gtilde_\rho, \mathfrak o)}
       {t-1}
=
  -{\mathbb T}^{K}_{\lambda}(\rho).
\]

Summarizing the above calculation,
we have shown that 
the rational function $\T(M_K, \gtilde_{\rho}, \mathfrak o)$
 has a simple zero at $t=1$ and 
its differential coefficient at $t=1$ agrees with 
minus the twisted Reidemeister torsion $-\mathbb T^{K}_{\lambda}(\rho)$.
\end{proof}

\subsection{
A description of $\mathbb T^{K}_{\lambda}$ 
using a Wirtinger representation}

Let $K$ be a knot in $S^3$
and $E_K$ its exterior.
We assume that $\rho \in R(\pi_1(E_K), G)$ is $\lambda$-regular. 
From Theorem \ref{thm:main_theorem}
we can describe $-\mathbb T^{K}_{\lambda}(\rho)$
by using the differential coefficient of $\T(E_K, \gtilde_{\rho}, \mathfrak o)$.
We will describe 
the differential coefficient of $\T(E_K, \gtilde_{\rho}, \mathfrak o)$ 
more explicitly
by using a Wirtinger representation of $\pi_1(E_K)$.

\medskip

\noindent
For a Wirtinger representation:
\[
\pi_1(E_K)=\langle x_1, \ldots, x_k \,|\, r_1, \ldots, r_{k-1}\rangle,
\]
we obtain a $2$-dimensional CW-complex $W$ which consists of 
one $0$-cell $p$, $k$ $1$-cells $x_1, \ldots, x_k$ and 
$(k-1)$ $2$-cells $D_1, \ldots, D_{k-1}$ attached by 
the relation $r_1, \ldots, r_{k-1}$.
This CW-complex $W$ is simple homotopy equivalent to $E_K$.
Let $\alpha : \pi_1(E_K) \to \Z=\langle t \rangle$ such that $\alpha(\mu)=t$.
Here $\mu$ is a meridian of $K$.
Note that
for all $i$, $\alpha(x_i)$ is equal to $t$ in $\Z=\langle t \rangle$.

The following calculation is due to the result of \cite{KL, Kitano}.
This chain complex $C_*(W; \gtilde_{\rho})$ is as follows:
\[
0 \to \mathfrakg(t)^{k-1} 
  \xrightarrow{\partial_2} \mathfrakg(t)^{k} 
  \xrightarrow{\partial_1} \mathfrakg(t) 
  \to 0 
\]
where
\begin{align*}
\partial_2
&=
\left(
\begin{array}{ccc}
\Phi(\frac{\partial r_1}{\partial x_1}) &\ldots & \Phi(\frac{\partial r_{k-1}}{\partial x_1}) \\
\vdots & \ddots & \vdots \\
\Phi(\frac{\partial r_1}{\partial x_k}) &\ldots & \Phi(\frac{\partial r_{k-1}}{\partial x_k}) 
\end{array}
\right),\\
& \\
\partial_1
&=
\left(
\Phi(x_1 -1),\, \Phi(x_2 -1),\, \ldots ,\, \Phi(x_k-1)
\right).
\end{align*}
Here we briefly denote the $l$-times direct sum of $\mathfrakg(t)$ 
by $\mathfrakg(t)^l$. 

We denote by $A^1_{K, Ad\circ \rho}$ $3(k-1)\times 3(k-1)$ matrix:
\[
\left(
\begin{array}{ccc}
\Phi(\frac{\partial r_1}{\partial x_2}) & \ldots & \Phi(\frac{\partial r_{k-1}}{\partial x_2}) \\
\vdots & \ddots & \vdots \\
\Phi(\frac{\partial r_1}{\partial x_k}) & \ldots & \Phi(\frac{\partial r_{k-1}}{\partial x_k}) 
\end{array}
\right).
\]
Under this situation, 
the twisted Alexander invariant $\T(W, \gtilde_{\rho}, \mathfrak o)$ is given by
\[
\tau_0 \cdot 
\frac{\det \, A^1_{K, Ad\circ \rho}}{\det (\Phi(x_1-1))}
\]
up to a factor $t^m\, (m \in \Z)$.

If $\rho(x_i)$ is conjugate to the upper triangulate matrix
\[
\left(
\begin{array}{cc}
 a & * \\
 0 & a^{-1} 
\end{array}
\right),
\]
then $Ad_{\rho(x_i^{-1})}$ is conjugate to the upper triangulate matrix
\[
\left(
\begin{array}{ccc}
1 &  *    & *  \\
  &  a^2  & *  \\
  &       & a^{-2}
\end{array}
\right).
\]
Calculating $\det (\Phi(x_1-1))$, we have that
\[
\det (\Phi(x_1-1))
=
(t-1)(t^2 - \trace(\rho(x^2_1))t +1).
\]
Since $\T(E_K, \gtilde_{\rho}, \mathfrak o)$ has zero at $t=1$, 
\begin{align*}
\left.
\frac{d}{dt}\T(E_K, \gtilde_{\rho}, \mathfrak o)
\right|_{t=1}
&=
\lim_{t \to 1}
  \frac{\T(E_K, \gtilde_{\rho}, \mathfrak o)}{t-1}\\
&=
\lim_{t \to 1}
  \tau_0\cdot t^m 
  \frac{\det A^1_{K, Ad\circ \rho}(t)}{(t-1)^2(t^2 - \trace(\rho(x_1^2))t +1)}.
\end{align*}

\begin{lemma}
If $\trace \rho(\partial E_K) \not \subset \{\pm 2\}$, then we have
\[
\lim_{t \to 1}\tau_0 \cdot t^m \frac{\det A^1_{K, Ad\circ \rho}(t)}{(t-1)^2}
=
\frac{\tau_0}{2} 
\left. \frac{d^2}{dt^2} \det A^1_{K, Ad\circ \rho}(t) \right|_{t=1}.
\]
\end{lemma}

\begin{proof}
The function $\T(E_K, \gtilde_{\rho}, \mathfrak o)$ has a simple zero at $t=1$ 
and 
the numerator $\det A^1_{K, Ad\circ \rho}(t)$ is an element of $\F[t, t^{-1}]$.
Hence 
$(t-1)^2$ divides $\det A^1_{K, Ad\circ \rho}(t)$.
We write $(t-1)^2 f(t)$ for $\det A^1_{K, Ad\circ \rho}(t)$.
Then the left hand side turns into $\lim_{t \to 1} \tau_0 \cdot t^m f(t)$, 
i.e., $\tau_0 f(1)$.
On the other hand, the right hand side becomes as follows.
\begin{align*}
\frac{\tau_0}{2} 
\left. \frac{d^2}{dt^2} \det A^1_{K, Ad\circ \rho}(t) \right|_{t=1}
&=
\frac{\tau_0}{2} \left. \frac{d^2}{dt^2} (t-1)^2 f(t) \right|_{t=1}\\
&=
\left.
  \frac{\tau_0}{2} 
  \frac{d}{dt} 
  \left\{
    2(t-1)f(t) 
    + (t-1)^2 f'(t)
  \right\} 
\right|_{t=1}\\
&=
\frac{\tau_0}{2} \left[ 2f(t) + 4(t-1)f'(t) + (t-1)^2 f''(t)\right]_{t=1} \\
&=\tau_0 f(1).
\end{align*}
\end{proof}

The numerator $\det A^1_{K, Ad\circ \rho}(t)$ is called 
{\it the first homology torsion\/} of $C_*(E_K; \gtilde_\rho)$ \cite{KL}.
We denote the first homology torsion by $\Delta_1(t)$.
By the above calculations, 
we obtain the following description of $\mathbb T^{K}_{\lambda}(\rho)$.
\begin{proposition}\label{another_form_torsion}
If $\trace(\rho (\partial E_K)) \not \subset \{\pm 2\}$, then we have
the following expression.
\[
\mathbb T^{K}_{\lambda}(\rho)
=
-
  \left.
  \frac{d}{dt}\T(E_K, \gtilde_\rho, \mathfrak o)
  \right|_{t=1}
=
  \frac{\tau_0 \Delta''_1(1) }{2}\cdot
  \frac{1}{\trace(\rho(x_1^2)) -2}.
\]
\end{proposition}

\begin{remark}
If $G$ is $\SU(2)$ and $\rho$ is $\lambda$-regular, 
then $\trace (\rho(\partial E_K)) \not \subset \{\pm 2\}$.
\end{remark}

\begin{remark}
We use a Wirtinger representation of $\pi_1(E_K)$ 
to describe $\T(E_K, \gtilde_{\rho}, \mathfrak o)$
in the above calculation.
The twisted Alexander invariant $\T(E_K, \gtilde_{\rho}, \mathfrak o)$
does not depend on 
the representation of $\pi_1(E_K)$ \cite{Wada}.
Since $\T(E_K, \gtilde_{\rho}, \mathfrak o)$
is determined by the finite presentable group $\pi_1(E_K)$ 
and $\rho \in R(E_K, G)$,
we do not necessarily need to use a Wirtinger representation
on calculating $\T(E_K, \gtilde_{\rho}, \mathfrak o)$.
\end{remark}

\section{Applications.}
\label{applications}
In this section, we deal with a $2$-bridge knot $K$ in $S^3$ 
and $\SU(2)$-representations of its knot group.
In this case
$\rho \in R(\pi_1(E_K), \SU(2))$ is irreducible 
if and only if 
$\rho(\pi_1(E_K))$ is a non-abelian subgroup of $\SU(2)$.
We will show 
the explicit calculation of $\SU(2)$-twisted Reidemeister torsion 
associated to $5_2$ knot
and study the critical points of the twisted Reidemeister torsion 
${\mathbb T}^K_{\lambda}$.
If $K$ is hyperbolic and $G$ is $\SL(2, \C)$,
then some features of ${\mathbb T}^K_{\mu}(\rho)$,
given in this section,
have appeared
in \cite[Section 4.3]{Porti}.

\subsection{A review of a representation of 
a $2$-bridge knot group}
It is well known that $\pi_1(E_K)$ has the representation:
\[
\langle 
  x, y 
  \,|\,
  wx=yw
\rangle,
\]
where $w$ is a word in $x$ and $y$.
Here $x$ and $y$ represent the meridian of the knot.
The method we use to describe the space of $\SL(2, \C)$ and 
$\SU(2)$-representations is due to 
R. Riley\,(\cite{Riley}). 
He shows how to parametrize conjugacy classes of irreducible $\SL(2, \C)$ and 
$\SU(2)$-representations of any $2$-bridge knot group.
We review his method (\cite{KK, Riley}).

\medskip

Given $s, u \in \C$, we consider the assignment as follows:
\[
x \mapsto 
\left(
\begin{array}{cc}
 s & 1 \\
 0 & 1
\end{array}
\right), 
\quad
y \mapsto 
\left(
\begin{array}{cc}
   s & 0 \\
 -su & 1
\end{array}
\right).
\] 
Let $W$ be the matrix obtained by replacing $x$ and $y$ 
by the above two matrices in the word $w$.
This assignment defines a $\GL(2, \C)$-representation 
if and only if $\phi(s, u)=0$ where $\phi(s, u)=W_{11} + (1-s)W_{12}$. 

One can obtain an $\SL(2, \C)$-representation 
from this $\GL(2, \C)$-representation by dividing the above two matrices 
by some square root of $s$.
If we give a path $(s(a), u(a))$ in $\C^2$ with $\phi(s(a), u(a))=0$ and 
some continuous branch of the square root along $s(a)$, 
then we obtain a path of $\SL(2, \C)$-representations.
Furthermore, 
all conjugacy classes of non-abelian $\SL(2, \C)$-representations arise 
in this way.

According to Proposition $4$ of Riley's paper \cite{Riley}, 
a pair $(s, u)$ with $\phi(s, u)=0$ corresponds to an $\SU(2)$-representation
if and only if 
$|s|=1$, and $u$ is real number which lies 
in the interval $[s+s^{-1}-2, 0] = [2\cos \theta -2, 0]$ 
where $s= e^{i \theta}$.
This correspondence means that 
the $\SL(2, \C)$-representation resulting from such a pair $(s, u)$ and 
some square root of $s$
is conjugate to an $\SU(2)$-representation in $\SL(2, \C)$.

\medskip

We take the ordered basis $E, H, F$ of ${\mathfrak{sl}}(2, \C)$ as follows.
\[
E
=
\left(
\begin{array}{cc}
0 & 1 \\
0 & 0
\end{array}
\right), 
H
=
\left(
\begin{array}{cc}
1 & \hphantom{-}0 \\
0 & -1
\end{array}
\right), 
F
=
\left(
\begin{array}{cc}
0 & 0 \\
1 & 0
\end{array}
\right).
\] 
The Lie algebra ${\mathfrak{su}}(2)$ is a subspace of ${\mathfrak{sl}}(2, \C)$.
The vectors $E, H, F$ also form a basis of ${\mathfrak{su}}(2)$.
Since the Euler characteristic of $E_K$ is zero, 
the non-abelian Reidemeister torsion $\mathbb T^K_{\lambda}(\rho)$ 
does not depend on a choice of a basis of 
${\mathfrak{su}}(2)$.
We can use $E, H, F$ as an ordered basis of ${\mathfrak{su}}(2)$.
We denote by $\rho_{\sqrt{s},u}$ the representation corresponding to the pair $(\sqrt{s}, u)$.
The representation matrices of $Ad (\rho_{\sqrt{s}, u}(x)) $ and 
$Ad(\rho_{\sqrt{s}, u}(y))$ 
for this ordered basis are given as follows.

\begin{lemma}\label{rep_matrix_adjoint}
\[
Ad(\rho_{\sqrt{s}, u}(x)) 
=
\left(
\begin{array}{ccc}
s & -2            & -\frac{1}{s} \\
0 & \hphantom{-}1 & \hphantom{-}\frac{1}{s} \\
0 & \hphantom{-}0 & \hphantom{-}\frac{1}{s}
\end{array}
\right), 
\quad 
Ad(\rho_{\sqrt{s}, u}(y))
=
\left(
\begin{array}{ccc}
s &  0  & 0  \\
su &  1 & 0 \\
-su^2 & -2u & \frac{1}{s}
\end{array}
\right).
\] 
\end{lemma}  

\noindent
Note that 
even if we choose another square root of $s$, 
we obtain the same representation matrices of 
the adjoint actions of $\rho_{\sqrt{s}, u}(x)$ and $\rho_{\sqrt{s}, u}(y)$.

\subsection{$\SU(2)$-twisted Reidemeister torsion associated to $5_2$ knot}
\begin{figure}[htb]
\begin{center}
\includegraphics{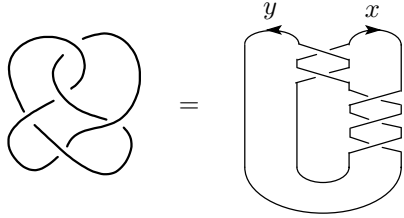}
\end{center}
\caption{A diagram of $5_2$ knot.} 
\label{diagram_5_2}
\end{figure}

We consider $5_2$ knot in the knot table of Rolfsen \cite{Rolfsen}.
Note that this knot is not fibered, since its Alexander polynomial is not monic.
This is the simplest example such as non-fibered in $2$-bridge knots.
Let $K$ be $5_2$ knot.
A diagram of $K$ is shown as in Figure \ref{diagram_5_2}.
%
%

This knot is also called $3$-twist knot.
It follows from Theorem $3$ of \cite{Klassen} 
that $\hat R^{irr}(\pi_1(E_K), \SU(2))$ consists of one circle and one open arc. 

\begin{figure}[htb]
\begin{center}
\includegraphics{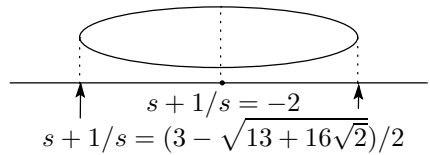}
\end{center}
\caption{$\hat R^{irr}(\pi_1(E_K), \SU(2))$ where $K$ is $5_2$ knot.} 
\label{fig:SU2_rep_5_2}
\end{figure}

The knot group $\pi_1 (E_K)$ has the following representation:
\[
\langle 
  x, y 
  \,|\,
  wx=yw
\rangle
\]
where 
$
  w=x^{-1}y^{-1}xyx^{-1}y^{-1}.
$
From this representation, the Riley's polynomial of $5_2$ is given by
\begin{align*}
& W_{11}+(1-s)W_{12} \\
&=
\frac{
-u^3 
+( 2( s+1/s )-3 )u^2 
+(-(s^2+1/s^2 )+3(s+1/s )-6 )u 
+2(s+1/s)-3
}
{s}.
\end{align*}
We may take Riley's polynomial $\phi(s, u)$ as 
\[
u^3 
-( 2( s+1/s)-3 )u^2 
+((s^2+1/s^2 )-3(s+1/s )+6 )u 
-(2(s+1/s)-3).
\]
We want to know pairs $(s, u)$ 
such that $s=e^{i \theta}$, $u$ is a real number 
in the interval $[2\cos \theta -2, 0]$
and $\phi(s, u)=0$.
When we regard $\phi(s, u)=0$ as the equation of $u$, 
the relation between the number of solutions of $\phi(s, u)=0$ and 
$s$ is as follows.
\begin{itemize}
\item[(1)] 
	If $-2 \leq s+1/s < (3-\sqrt{13+16\sqrt{2}})/2$, 
	then $\phi(s, u)=0$ has three different simple root in $[s+1/s -2, 0]$.
\item[(2)] 
	If $s+1/s = (3-\sqrt{13+16\sqrt{2}})/2$, 
	then $\phi(s, u)=0$ has a simple root and a multiple root in 
        $[s+1/s -2, 0]$.
\item[(3)]
	If $(3-\sqrt{13+16\sqrt{2}})/2 <  s+1/s < 3/2$,
	then $\phi(s, u)=0$ has a simple root in $[s+1/s -2, 0]$.
\end{itemize}
The figure of $\hat R^{irr}(\pi_1(E_K), \SU(2))$ is given as in 
Figure \ref{fig:SU2_rep_5_2}.
%
%

We denote the $\SU(2)$-representation corresponding to $(s, u)$ 
by $\rho_{\sqrt{s}, u}$.
Then we can express $\mathbb T^K_{\lambda}(\rho_{\sqrt{s}, u})$ 
from Proposition \ref{another_form_torsion} 
as follows.
\[
\mathbb T^K_{\lambda}(\rho_{\sqrt{s}, u})
=
\frac{\tau_0 \Delta''_1(1)}{2} \cdot \frac{1}{s+1/s -2}
\]

Using a computer,
we calculate 
a half of the differential coefficient of the second order of the numerator
and simplify with the equation $\phi(s, u)=0$.
Then we have
\begin{align*}
& \frac{\tau_0 \Delta''_1(1)}{2}\\
&=\tau_0
(s+1/s -2)
(-(5(s+1/s)+3)u^2 +(5(s+1/s)^2-7(s+1/s)+1)u+1-10(s+1/s)).
\end{align*}
Therefore we have
\[
\mathbb T^K_{\lambda}(\rho_{\sqrt{s}, u})
=
\tau_0
(-(5(s+1/s)+3)u^2 +(5(s+1/s)^2-7(s+1/s)+1)u+1-10(s+1/s)),
\]
where $(u, s)$ satisfies $\phi(u, s)=0$.

\subsection{
On critical points of the $\SU(2)$-twisted Reidemeister torsion 
associated to $2$-bridge knots}
\label{torsion_2_bridge}
From the example in the previous subsection, 
one can guess that 
the $\SU(2)$-twisted Reidemeister torsion 
$\mathbb T^{K}_{\lambda}$
associated to a $2$-bridge knot $K$ is 
a function for the parameter $s + 1/s$.
Indeed the following holds.

\begin{proposition}\label{prop:parameter_torsion}
Let $K$ be a $2$-bridge knot and $\gamma$ a simple closed curve in
the boundary torus of $E_K$.
Suppose that $\gamma$-regular $\SU(2)$-representations are parametrized 
by $(s, u) \in \U(1) \times \R$ of Riley's method.
If the trace of the meridian, $\sqrt{s} + 1/\sqrt{s}$, gives a local parameter of the $\SU(2)$-character variety,
then the twisted Reidemeister torsion $\mathbb T^{K}_{\gamma}$ 
is a smooth function for $s+1/s$.
\end{proposition}  

\begin{proof}
If we denote by $\rho_{\sqrt{s}, u}$ a $\gamma$-regular representation
corresponding to $\sqrt{s}+1/\sqrt{s}$, then there exists some homomorphism
$\varepsilon: \pi_1(E_K) \to \{\pm 1\}$ such that 
$\varepsilon \rho_{\sqrt{s}, u}$ is a $\gamma$-regular representation 
corresponding to $-\sqrt{s}-1/\sqrt{s}$.
By the construction of ${\mathbb T}^K_{\gamma}$, 
${\mathbb T}^K_{\gamma}(\rho)$ 
is equal to 
${\mathbb T}^K_{\gamma}(\varepsilon \rho)$.
Since $\sqrt{s} + 1/\sqrt{s}$ is a square root of $s + 1/s + 2$ and regular
representations are irreducible,
the twisted Reidemeister torsion ${\mathbb T}^{K}_{\gamma}$ 
is a smooth function for $s + 1/s$.
\end{proof}

\begin{cor}\label{cor:parameter_torsion}
If the trace of the meridian gives a local parameter of the $\SU(2)$-character variety and the twisted Reidemeister torsion ${\mathbb T}^K_{\lambda}$ is defined,
then ${\mathbb T}^K_{\lambda}$ is a smooth function for $s+1/s$.
\end{cor}

\begin{remark}\label{rem:parameter_dihedral}
All representations $\rho$ of $2$-bridge knot groups into $\SU(2)$
such that $\trace(\rho(\mu)) = 0$ 
are binary dihedral representations.
It follows from \cite{HK} that there exists a neighbourhood of the character of each binary dihedral representation for any $2$-bridge knot, which is diffeomorphic to an open interval.
From \cite{Burde}, the trace of the meridian gives a local parameter on a neighbourhood of the character of each dihedral representation for $2$-bridge knots.
\end{remark}

We can regard the twisted Reidemeister torsion $\mathbb T^K_{\lambda}$ as a smooth function on a neighbourhood of the character of each binary dihedral representation.
Moreover these characters can be critical points of $\mathbb T^K_{\lambda}$ as follows.

\begin{cor}\label{critical_pt_torsion}
Let $K$ be a $2$-bridge knot. 
If a $\lambda$-regular component of 
the $\SU(2)$-character variety of $\pi_1(E_K)$
contains the characters of dihedral representations, 
then the function $\mathbb T^K_{\lambda}$ has a critical point at the character of each dihedral representation.
\end{cor}

\begin{proof}
By Corollary \ref{cor:parameter_torsion} and Remark \ref{rem:parameter_dihedral},  
the twisted Reidemeister torsion $\mathbb T^K_{\lambda}$ is a smooth function for $s + 1/s$.
When we substitute $e^{i \theta}$ for $s$, 
we can describe $\mathbb T^K_{\lambda}(\rho)$  as
\[
  \frac{f(2 \cos \theta)}{2\cos \theta -2}
\]
where $f(2 \cos \theta)$ is a smooth function for $2 \cos \theta$.
This is a description of $\mathbb T^K_{\lambda}$ 
with respect to the local coordinate $\theta$
of $\hat R^{irr}(\pi_1(E_K), \SU(2))$.
The derivation of this function for $\theta$ becomes 
\[
\frac{\{-2 f'(2 \cos \theta)(2\cos \theta -2)+2f(2 \cos \theta) \} \sin \theta}
     {(2\cos\theta -2)^2}.
\]
We recall that 
$\trace(\rho_{\sqrt{s}, u}(\mu)) = \trace(\rho_{\sqrt{s}, u}(x))=2\cos (\theta/2)$.
If $\trace(\rho_{\sqrt{s}, u}(\mu))=2\cos (\theta/2) = 0$, then $\sin \theta =0$.
Hence the derivation of $\mathbb T^K_{\lambda}$ vanishes 
if $\rho$ satisfies $\trace(\rho(\mu))=0$.
\end{proof}

\begin{remark}
From \cite{Burde}, for $2$-bridge knots,
the character of a binary dihedral representation is a branch point 
of the two-fold branched cover from the $\SU(2)$-character variety
to the ${\rm SO}(3)$-character variety.
Moreover, every algebraic component of the $\SU(2)$-character variety
contains the character of such a representation.
\end{remark}

\begin{remark}
By \cite[Theorem 10]{Klassen}, for a knot $K$,
the number of conjugacy class of binary dihedral representations is given by
$(|\Delta_K(-1)|-1)/2$
where $\Delta_K(t)$ is the Alexander polynomial of $K$.
In particular, for a $2$-bridge knot $b(\alpha, \beta)$ 
(Schubert's notation, see for example \cite{Knots}),
this number is given by $(\alpha -1)/2$.
\end{remark}


\section*{Acknowledgements}
The author would like to express sincere gratitude to Mikio Furuta 
for his suggestions and helpful discussions.
He is thankful to Hiroshi Goda, Takayuki Morifuji, Teruaki Kitano 
for helpful suggestions.
Especially the author gratefully acknowledges the many helpful suggestions of
Hiroshi Goda during the preparation of the paper.
Our main theorem was written as the statement for knots in $S^3$ at first. 
J\'er\^ome Dubois pointed out that 
our main theorem can hold for knots in homology three spheres.
The author is thankful to J\'er\^ome Dubois for his pointing out.
He also would like to thank the referee for his/her careful reading 
and appropriate advices.
He/She has given suggestions to improve the proofs of 
Proposition \ref{PropositionA} and Proposition \ref{prop:parameter_torsion}.
He/She also suggested the fact that critical points of Reidemeister torsion 
are related to the dihedral representations.


\end{document}